\documentclass[12pt]{article}
\usepackage{amssymb}
\usepackage{graphicx}
\usepackage{amsmath}

\newtheorem{theorem}{Theorem}

\newtheorem{lemma}[theorem]{Lemma}

\newtheorem{proposition}[theorem]{Proposition}

\begin{document}

\title{On the number of generators needed for free profinite products of finite
groups}
\author{Mikl\'{o}s Ab\'{e}rt and P\'{a}l Heged\H us
\footnote{The first author's research is partially supported by NSF grant DMS-0401006. The second author's research is partially supported by OTKA grant T38059 and the Magyary Zolt\'{a}n Postdoctoral Fellowship}
}
\maketitle

\begin{abstract}
We provide lower estimates on the minimal number of generators of the
profinite completion of free products of finite groups.

In particular, we show that if $C_{1},\ldots ,C_{n}$ are finite cyclic
groups then there exists a finite group $G$ which is generated by isomorphic
copies of $C_{1},\ldots ,C_{n}$ and the minimal number of generators of $G$
is $n$.
\end{abstract}

\section{Introduction}

For a group $G$ let $d(G)$ denote the minimal number of generators for $G$.
If $G$ is a profinite group then we mean topological generation rather than
the abstract one. Let $\widehat{G}$ denote the profinite completion of $G$;
trivially $d(\widehat{G})\leq d(G)$. The first finitely generated residually
finite examples where the two quantities are different were found by Noskov 
\cite{noskov}. His examples were metabelian and he also showed that for
these groups we have 
\begin{equation*}
d(G)\leq (t^{2}+5t+2)/2\text{ with }t=d(\widehat{G}).
\end{equation*}
An old question of Melnikov \cite[6.31]{melnik} asked whether $d(G)$ is
always bounded by a function of $d(\widehat{G})$ for a residually finite,
finitely generated group $G$. This has been recently answered negatively by
Wise \cite{wise} but is still open for linear groups.

Another class of groups where passing to profinite completion may imply a
drop in the minimal number of generators is free products of finite groups.
Let $G_{1},G_{2},\ldots ,G_{n}$ be finite groups, let $s=\max_{i}d\left(
G_{i}\right) $ and let 
\begin{equation*}
\Gamma =G_{1}\ast G_{2}\ast \cdots \ast G_{n}
\end{equation*}
be the free product of the $G_{i}$. Then the so-called Grushko-Neumann
theorem (see \cite{gru} and \cite{neu}) says that $d\left( \Gamma \right)
=\sum_{i}d(G_{i})$. On the other hand, Kov\'{a}cs and Sim \cite{kovsim}
showed that if the $G_{i}$ are solvable and have pairwise coprime orders
then $d\left( \widetilde{\Gamma }\right) \leq n+s-1$, where $\widetilde{%
\Gamma }$ denotes the prosolvable completion of $\Gamma $. In the languange
of finite groups this translates to stating that if $G$ is a finite solvable
group which is generated by subgroups isomorphic to $G_{1},G_{2},\ldots
,G_{n}$, then $d(G)\leq n+s-1$ (see \cite{riwong}).

This was followed by work of Lucchini \cite{luch1} who, using the
Classification of Finite Simple Groups, showed that there exists an absolute
constant $c$ such that if the $G_{i}$ have pairwise coprime orders and $n>2$
then 
\begin{equation*}
d\left( \widehat{\Gamma }\right) \leq (1+\frac{4c}{3})(n-1)+2s+c.
\end{equation*}
It is conjectured that if the $G_{i}$ have pairwise coprime orders then in
fact 
\begin{equation*}
d\left( \widehat{\Gamma }\right) =n+s-1.
\end{equation*}
The upper bound $d\left( \widehat{\Gamma }\right) \leq n+s$ is proved by
Lucchini in \cite[Theorem C]{luch1} in the case when the $G_{i}$ are $p_{i}$%
-groups for distinct primes $p_{i}$.

The aim of this paper is to support general lower bounds for $d(\widehat{%
\Gamma })$. For some special families of finite groups this has been done by
Kovacs and Sim \cite{kovsim}. The first estimate of this type works for
arbitrary finite groups and trivially implies $d\left( \Gamma \right) \leq d(%
\widehat{\Gamma })^{2}$.

\begin{theorem}
\label{en}Let $G_{1},G_{2},\ldots ,G_{n}$ be finite groups and let $\Gamma
=G_{1}\ast G_{2}\ast \cdots \ast G_{n}$. Then 
\begin{equation*}
d\left( \widehat{\Gamma }\right) \geq n.
\end{equation*}
\end{theorem}

In particular, if all the $G_{i}$ are nontrivial cyclic, then we have the
equality $d\left( \widehat{\Gamma }\right) =n$, proving the above conjecture
for the case $s=1$. In fact, as Proposition \ref{solsol} shows, already $%
d\left( \widetilde{\Gamma }\right) =n$.

Note that the weaker estimate

\begin{equation*}
d\left( \widehat{\Gamma }\right) \geq n-\sum_{i=1}^{n}\frac{1}{\left|
G_{i}\right| }
\end{equation*}
can be proved in various ways. E.g. it is immediate from the following
observation, which may be interesting in itself.

\begin{proposition}
\label{betti}Let $\Gamma $ be a finitely presented residually finite group.
Then 
\begin{equation*}
d\left( \widehat{\Gamma }\right) \geq b_{1}^{(2)}(\Gamma )+1
\end{equation*}
where $b_{1}^{(2)}(\Gamma )$ denotes the first $L_{2}$-Betti number of $%
\Gamma $.
\end{proposition}

Note that by the definition of $L_{2}$-Betti numbers $d\left( \Gamma \right)
\geq b_{1}^{(2)}(\Gamma )+1$ for arbitrary groups \cite{luck}.

Our second theorem involves $s$ in the lower estimate in the following form.

\begin{theorem}
\label{bonyi}Let $G_{1},G_{2},\ldots ,G_{n}$ be finite groups, let $\Gamma
=G_{1}\ast G_{2}\ast \cdots \ast G_{n}$ and let 
\begin{equation*}
s^{\prime }=\max (d(G_{i}/G_{i}^{\prime }))
\end{equation*}
Then 
\begin{equation*}
d\left( \widehat{\Gamma }\right) \geq n+s^{\prime }-1.
\end{equation*}
\end{theorem}

In particular, if all the $G_{i}$ are nilpotent then $d\left( \widehat{%
\Gamma }\right) \geq n+s-1$ which sets the conjectured lower bound. If
moreover the $G_{i}$ are $p_{i}$-groups for distinct primes $p_{i}$, then
using Lucchini's upper bound we get 
\begin{equation*}
n+s-1\leq d\left( \widehat{\Gamma }\right) \leq n+s.
\end{equation*}

For groups of pairwise coprime order where the minimal number of generators
is not witnessed by the abelianization, we are unable to set the conjectured
lower bound in general, but in the case $s=2$ we can show that it is the
best possible one can hope for.

\begin{theorem}
\label{pelda}For every $n$ there exist solvable groups $G_{1},G_{2},\ldots
,G_{n}$ of pairwise coprime order such that $d(G_{i})=2$, $%
G_{i}/G_{i}^{\prime }$ is cyclic ($1\leq i\leq n$) and for $\Gamma
=G_{1}\ast G_{2}\ast \cdots \ast G_{n}$ we have $d\left( \widehat{\Gamma }%
\right) \geq n+1$.
\end{theorem}

\section{Proofs}

First we prove Proposition \ref{betti}. Note that this is independent of the
rest of the paper and it provides a weaker bound than the one obtained with
our main method. However, it readily generalizes to all classes of groups
where we can compute the first $L_{2}$ Betti number, e.g., to amalgamated
products.

\bigskip

\noindent \textbf{Proof of Proposition \ref{betti}. }Let $%
N_{1}\vartriangleleft \Gamma $ be a normal subgroup of finite index such
that $d\left( \widehat{\Gamma }\right) =d(\Gamma /N)$. Let 
\begin{equation*}
\Gamma =N_{0}\geq N_{1}\geq N_{2}\geq \ldots
\end{equation*}
be an infinite chain of normal subgroups of $\Gamma $ of finite index such
that $\cap _{i}N_{i}=1$. Let $K_{i}=N_{i}^{\prime }N_{i}^{2}$, where $%
N_{i}^{\prime }$ denotes the derived subgroup and $N_{i}^{2}$ is the normal
subgroup generated by all squares in $N_{i}$ ($i\geq 0$). Let $G_{i}=\Gamma
/N_{i}$ and let $H_{i}=\Gamma /K_{i}$ ($i\geq 0$). Let $d_{i}$ denote the
torsion-free rank of the abelianization of $N_{i}$ (or in other words, the
first homology of $N_{i}$). Using a theorem of L\"{u}ck \cite{luck}, we have 
\begin{equation*}
\beta _{1}^{2}(\Gamma )=\lim_{n\rightarrow \infty }\frac{d_{n}}{\left|
G_{n}\right| }
\end{equation*}

Now $N_{i}/K_{i}$ is an elementary Abelian $2$-group and 
\begin{equation*}
d(N_{i}/K_{i})\geq d_{i}\text{ (}i\geq 0\text{)}.
\end{equation*}
The index of $N_{i}/K_{i}$ in $H_{i}$ is $\left| G_{i}\right| $ and so using
the Nielsen-Schreier theorem we have 
\begin{equation*}
d(N_{i}/K_{i})\leq (d(H_{i})-1)\left| G_{i}\right| +1\text{ (}i\geq 0\text{)}
\end{equation*}
which gives us 
\begin{equation*}
d(H_{i})\geq \frac{d(N_{i}/K_{i})-1}{\left| G_{i}\right| }+1\geq \frac{%
d_{i}-1}{\left| G_{i}\right| }+1\text{ (}i\geq 0\text{)}
\end{equation*}
Since $\cap _{i}K_{i}\leq \cap _{i}N_{i}=1$ and $\left| G_{n}\right|
\rightarrow \infty $, we have 
\begin{equation*}
d\left( \widehat{\Gamma }\right) =\lim_{n\rightarrow \infty }d(H_{i})\geq
\lim_{n\rightarrow \infty }\frac{d_{n}-1}{\left| G_{n}\right| }+1=\beta
_{1}^{2}(\Gamma )+1
\end{equation*}
The proposition holds. $\square $

\bigskip

Now we start building towards Theorem \ref{en} and Theorem \ref{bonyi}.

\bigskip

Let $\Gamma $ be a finitely generated group and $H$ a finite group. Let $%
\mathrm{Hom}(\Gamma ,H)$ denote the set of homomorphisms from $\Gamma $ to $%
H $. Then $\mathrm{Hom}(\Gamma ,H)$ is finite. Let 
\begin{equation*}
h(\Gamma ,H)=\frac{\log \left| \mathrm{Hom}(\Gamma ,H)\right| }{\log \left|
H\right| }
\end{equation*}
The number $h(\Gamma ,H)$ will be the key notion of this paper. Let 
\begin{equation*}
K(\Gamma ,H)=\bigcap_{\varphi \in \mathrm{Hom}(\Gamma ,H)}\ker \varphi
\end{equation*}
and let the quotient group 
\begin{equation*}
G(\Gamma ,H)=\Gamma /K(\Gamma ,H).
\end{equation*}

Since $K(\Gamma ,H)$ can be obtained as a finite intersection of subgroups
of finite index, $G(\Gamma ,H)$ is a finite image of $\Gamma $. Also, each
homomorphism from $\Gamma $ to $H$ factors through $K(\Gamma ,H)$, so we
have 
\begin{equation*}
\left| \mathrm{Hom}(\Gamma ,H)\right| =\left| \mathrm{Hom}(G(\Gamma
,H),H)\right|
\end{equation*}
implying 
\begin{equation*}
h(G(\Gamma ,H),H)=h(\Gamma ,H).
\end{equation*}

The following two basic lemmas are needed later.

\begin{lemma}
\label{free}Let $\Gamma _{i}$ ($1\leq i\leq n$) be finitely generated groups
and let $H$ be a finite group. Then 
\begin{equation*}
h(\Gamma _{1}\ast \Gamma _{2}\ast \cdots \ast \Gamma
_{n},H)=\sum_{i=1}^{n}h(\Gamma _{i},H)
\end{equation*}
\end{lemma}

\noindent \textbf{Proof. }By the definition of a free product, for every set
of homomorphisms $\varphi _{i}\in \mathrm{Hom}(\Gamma _{i},H)$ ($1\leq i\leq
n$) there exists a unique homomorphism $\varphi \in \mathrm{Hom}(\Gamma
_{1}\ast \Gamma _{2}\ast \cdots \ast \Gamma _{n},H)$ such that the
restriction of $\varphi $ to $\Gamma _{i}$ equals $\varphi _{i}$ ($1\leq
i\leq n$). Hence

\begin{equation*}
\left| \mathrm{Hom}(\Gamma _{1}\ast \Gamma _{2}\ast \cdots \ast \Gamma
_{n},H)\right| =\prod_{i=1}^{n}\left| \mathrm{Hom}(\Gamma _{i},H)\right|
\end{equation*}
implying 
\begin{eqnarray*}
h(\Gamma _{1}\ast \Gamma _{2}\ast \cdots \ast \Gamma _{n},H) &=&\frac{\log
\left| \mathrm{Hom}(\Gamma _{1}\ast \Gamma _{2}\ast \cdots \ast \Gamma
_{n},H)\right| }{\log \left| H\right| }= \\
&=&\frac{\sum_{i=1}^{n}\log (\left| \mathrm{Hom}(\Gamma _{i},H)\right| )}{%
\log \left| H\right| }=\sum_{i=1}^{n}h(\Gamma _{i},H)
\end{eqnarray*}
as claimed. $\square $

\begin{lemma}
Let $\Gamma $ be a finitely generated group and let $H$ be a finite group.
Then 
\begin{equation*}
h(\Gamma ,H^{n})=h(\Gamma ,H)
\end{equation*}
for all natural numbers $n$.
\end{lemma}

\noindent \textbf{Proof. }A function $\varphi :\Gamma \rightarrow H^{n}$ is
a homomorphism if and only if all the coordinate functions of $\varphi $ are
homomorphisms into $H$. Thus $\left| \mathrm{Hom}(\Gamma ,H^{n})\right|
=\left| \mathrm{Hom}(\Gamma ,H)\right| ^{n}$ which implies the statement. $%
\square $

\bigskip

The following lemma establishes a connection between the function $h$ and
the minimal number of generators for the profinite completion.

\begin{lemma}
\label{becs}Let $G_{i}$ be finite groups ($1\leq i\leq n$) and let 
\begin{equation*}
\Gamma =G_{1}\ast G_{2}\ast \cdots \ast G_{n}.
\end{equation*}
Then 
\begin{equation*}
d(\widehat{\Gamma })\geq \sum_{i=1}^{n}h(G_{i},H)
\end{equation*}
for any finite group $H$.
\end{lemma}

\noindent \textbf{Proof. }If $G$ is an arbitrary homomorphic image of $%
\Gamma $ then $\mathrm{Hom}(G,H)\leq \mathrm{Hom}(\Gamma ,H)$ and so we have 
$h(G,H)\leq h(\Gamma ,H)$ as well. In particular for $d=d(G)$ we have 
\begin{equation*}
h(G,H)\leq h(F_{d},H)=\frac{\log \left| \mathrm{Hom}(F_{d},H)\right| }{\log
\left| H\right| }=\frac{\log (\left| H\right| ^{d})}{\log \left| H\right| }=d
\end{equation*}

Using this and Lemma \ref{free} we have 
\begin{equation*}
d(\widehat{\Gamma })\geq d(G(\Gamma ,H))\geq h(G(\Gamma ,H),H)=h(\Gamma
,H)=\sum_{i=1}^{n}h(G_{i},H)
\end{equation*}
as claimed. $\square $

\bigskip

So in order to obtain a lower bound on $d(\widehat{\Gamma })$ we have to
find a target group $H$, such that all the $G_{i}$ have many homomorphisms
into $H$. Note that if we choose the target group to be a large symmetric
group or a large dimensional general linear group over a fixed finite field,
we get the estimate 
\begin{equation*}
d\left( \widehat{\Gamma }\right) \geq n-\sum_{i=1}^{n}\frac{1}{\left|
G_{i}\right| }
\end{equation*}
already established by Proposition \ref{betti}.\textbf{\ }It turns out that
the best target groups for our purposes will be produced from semisimple $%
G_{i}$-modules over finite fields.

\bigskip

\noindent \textbf{Proof of Theorem \ref{en}. }Recall that $O_{p}(G)$ denotes
the largest normal $p$-subgroup of $G$. Let $p$ be a prime such that 
\begin{equation*}
O_{p}(G_{i})=1\text{ }(1\leq i\leq n)
\end{equation*}
and let $F=\mathbb{F}_{p}$ be the field of order $p$. Let $M_{i}$ be a
nontrivial simple $G_{i}$-module over $F$ of dimension $d_{i}=\dim _{F}M_{i}$
($1\leq i\leq n$). Let $l$ be the least common multiple of the $d_{i}$ and
let $V$ be a vectorspace over $F$ of dimension $l$.

Let $1\leq i\leq n$. Since $d_{i}$ divides $l$, $V$ can be turned into a
semisimple $G_{i}$-module such that all the simple factors of $V$ under $%
G_{i}$ are isomorphic to $M_{i}$. Let $L_{i}\subseteq GL(V)$ denote the
linear action of $G_{i}$ on $V$. Since $M_{i}$ is nontrivial, $L_{i}$ is not
the trivial group and since $M_{i}$ is simple, $L_{i}$ has no nonzero fixed
vector in $V$.

Let 
\begin{equation*}
R=\left\langle L_{i}\mid 1\leq i\leq n\right\rangle \subseteq GL(V)
\end{equation*}
be the linear group generated by the $L_{i}$. Then $V$ is a semisimple $R$%
-module. Let $r=\left| R\right| $.

Let $m$ be a natural number and let $H$ be the semidirect product of $V^{m}$
and $R$. Then $H$ has order $p^{lm}r$. We want to estimate $\left| \mathrm{%
Hom}(G_{i},H)\right| $ from below. It will suffice to consider conjugates of
a fixed surjective homomorphism from $G_{i}$ to $L_{i}$. The number of those
conjugates equals the size of the conjugacy class of $L_{i}$ in $H$. Since $%
L_{i}$ has no fixed vector in $V$, the centralizer $Z_{H}(L_{i})\leq R$.
This implies 
\begin{equation*}
\left| \mathrm{Hom}(G_{i},H)\right| \geq \frac{\left| H\right| }{\left|
Z_{H}(L_{i})\right| }\geq \frac{\left| H\right| }{r}=p^{lm}
\end{equation*}
so 
\begin{equation*}
h(G_{i},H)\geq \frac{\log p^{lm}}{\log (p^{lm}r)}=1-\frac{\log r}{m\log
p^{l}+\log r}
\end{equation*}
Using Lemma \ref{becs} this gives 
\begin{equation*}
d(\widehat{\Gamma })\geq \sum_{i=1}^{n}h(G_{i},H)\geq n\left( 1-\frac{\log r%
}{m\log p^{l}+\log r}\right)
\end{equation*}
Letting $m$ to be arbitrarily large this leads to 
\begin{equation*}
d(\widehat{\Gamma })\geq n
\end{equation*}
The theorem holds. $\square $

\bigskip

Now we prove Theorem \ref{bonyi} using the construction above.

\bigskip

\noindent \textbf{Proof of Theorem \ref{bonyi}. }We can assume that $%
d(G_{n}/G_{n}^{\prime })=s^{\prime }$. Let $p$ be a prime such that $%
p^{s^{\prime }}$ divides $\left| G_{n}/G_{n}^{\prime }\right| $ and let $F=%
\mathbb{F}_{p}$ be the field of order $p$. By permuting the $G_{i}$ we can
also assume that there exists $0\leq t<n$ such that $p$ does not divide $%
\left| G_{i}/G_{i}^{\prime }\right| $ ($1\leq i\leq t$) and $p$ divides $%
\left| G_{i}/G_{i}^{\prime }\right| $ ($t+1\leq i\leq n$). Let us define a
new list of finite groups $H_{i}$ ($1\leq i\leq t+1$) as follows. For $1\leq
i\leq t$ let 
\begin{equation*}
H_{i}=G_{i}/O_{p}(G_{i})
\end{equation*}
and let 
\begin{equation*}
H_{t+1}=C_{p}^{s^{\prime }+n-t-1}\text{.}
\end{equation*}

From here we follow the construction and notation in the proof of Theorem 
\ref{en} using the $H_{i}$ ($1\leq i\leq t$) and $p$ as prime. This is
allowed since $O_{p}(H_{i})=1$ ($1\leq i\leq t$). For a large enough $m$ let 
$H$ be the target group given by the construction. Then, as before, we have 
\begin{equation*}
h(H_{i},H)\geq 1-\frac{\log r}{m\log p^{l}+\log r}\text{ }(1\leq i\leq t)%
\text{.}
\end{equation*}
Now $\mathrm{Hom}(H_{t+1},V^{m})\subseteq \mathrm{Hom}(H_{t+1},H)$ implying 
\begin{eqnarray*}
h(H_{t+1},H) &=&\frac{\log \left| \mathrm{Hom}(H_{t+1},H)\right| }{\log
\left| H\right| }\geq \\
&\geq &\frac{\log \left| \mathrm{Hom}(H_{t+1},V^{m})\right| }{m\log \left|
V\right| +\log \left| S\right| }=h(H_{t+1},V^{m})\frac{m\log \left| V\right| 
}{m\log \left| V\right| +\log \left| S\right| }
\end{eqnarray*}
and $\left| \mathrm{Hom}(H_{t+1},C_{p})\right| =p^{s^{\prime }+n-t-1}$
implying 
\begin{equation*}
h(H_{t+1},V^{m})=h(H_{t+1},C_{p})=s^{\prime }+n-t-1
\end{equation*}
which gives 
\begin{equation*}
\sum_{i=1}^{t+1}h(H_{i},H)\geq (s^{\prime }+n-t-1)\frac{m\log \left|
V\right| }{m\log \left| V\right| +\log \left| S\right| }+t\left( 1-\frac{%
\log r}{m\log p^{l}+\log r}\right)
\end{equation*}
Setting $m$ to be large enough and using Lemma \ref{becs} we get 
\begin{equation*}
d(\widehat{H_{1}\ast \cdots \ast H_{t+1}})\geq
\sum_{i=1}^{t+1}h(H_{i},H)\geq s^{\prime }+n-1.
\end{equation*}
On the other hand, $H_{i}$ is a quotient of $G_{i}$ ($1\leq i\leq t$) and $%
H_{t+1}$ is a quotient of $G_{t+1}\ast \cdots \ast G_{n}$ which implies that 
$H_{1}\ast \cdots \ast H_{t+1}$ is a quotient of $\Gamma $, leading to 
\begin{equation*}
d(\widehat{\Gamma })\geq d(\widehat{H_{1}\ast \cdots \ast H_{t+1}})\geq
s^{\prime }+n-1\text{. }
\end{equation*}
The theorem holds. $\square $

\bigskip

Now we prove Proposition \ref{solsol}. It is again a slight modification of
the construction in Theorem \ref{en}.

\begin{proposition}
\label{solsol}Let $G_{1},G_{2},\ldots ,G_{n}$ be finite cyclic groups and
let $\Gamma =G_{1}\ast G_{2}\ast \cdots \ast G_{n}$. Then 
\begin{equation*}
d\left( \widetilde{\Gamma }\right) =d\left( \widehat{\Gamma }\right) =n
\end{equation*}
where $\widetilde{\Gamma }$ denotes the prosolvable completion of $\Gamma $.
\end{proposition}

\noindent \textbf{Proof. }Obviously $d\left( \widetilde{\Gamma }\right) \leq
d\left( \widehat{\Gamma }\right) \leq n$ holds, so it is enough to show that 
$d\left( \widetilde{\Gamma }\right) \geq n$. Just as before, we can assume
that the $G_{i}$ have prime order $p_{i}$ ($1\leq i\leq n$). Let 
\begin{equation*}
k=\prod_{i=1}^{n}p_{i}\text{.}
\end{equation*}
By Dirichlet's theorem there are infinitely many primes in the arithmetic
progression $kn+1$ ($n\in \mathbb{N}$). Let $p$ be such a prime and let $%
F=F_{p}$. Then the multiplicative group $F^{\ast }$ is a cyclic group of
order divisible by all the $p_{i}$, so $G_{i}$ embeds nontrivially into $%
F^{\ast }$ ($1\leq i\leq n$). In other terms, $F$ can be turned into a
nontrivial one-dimensional $G_{i}$-module. The linear actions of the $G_{i}$
will generate a subgroup of $F^{\ast }$ that is isomorphic to the direct
product $G$ of the $G_{i}$. Following the construction, we get a target
group $H$ which is metabelian, being the extension of an $F$-space by $G$.
But then the witness group 
\begin{equation*}
G(\Gamma ,H)=\Gamma /K(\Gamma ,H)=\Gamma /\bigcap_{\varphi \in \mathrm{Hom}%
(\Gamma ,H)}\ker \varphi \hookrightarrow H^{\left| \mathrm{Hom}(\Gamma
,H)\right| }
\end{equation*}
embeds into the direct product of metabelian groups, so it is metabelian
itself. By Lemma \ref{becs} we have $d(G(\Gamma ,H))\geq n$, finishing the
proof. $\square $

\bigskip

Now we start building towards Theorem \ref{pelda}. The first result needed
is due to Erd\H os \cite{erdos} and is purely number-theoretic.

\begin{theorem}[Erd\H os]
\label{szamolos} Let $A$ be an infinite set of positive integers and let 
\begin{equation*}
f_{n}(A)=\left| A\cap \left\{ 1,\ldots ,n\right\} \right| \text{. }
\end{equation*}
Assume that \newline
1) $f_{n}(A)$ increases faster than $n^{(\surd 5-1)/2}$; \newline
2) Every arithmetic progression contains at least one integer which is the
sum of distinct elements of $A$. \newline
Then every sufficiently large integer is a sum of distinct elements of $A$.
\end{theorem}

For a prime $p$ let $\mathrm{Aff}(p)$ denote the group of affine
transformations of $\mathbb{F}_{p}$. Then $\mathrm{Aff}(p)$ acts on $\mathbb{%
F}_{p}$ and so it embeds into the symmetric group $\mathrm{Sym}(\mathbb{F}%
_{p})$.

\begin{lemma}
\label{pici}Let $H\leq \mathrm{Aff}(p)$ be a subgroup properly containing
the additive subgroup $\mathbb{F}_{p}$. Then the centralizer $Z_{\mathrm{Sym}%
(\mathbb{F}_{p})}(H)=1$.
\end{lemma}

\noindent \textbf{Proof. }Since $\mathbb{F}_{p}$ is Abelian and transitive
in $\mathrm{Sym}(\mathbb{F}_{p})$, $Z_{\mathrm{Sym}(\mathbb{F}_{p})}(\mathbb{%
F}_{p})=\mathbb{F}_{p}$, implying $Z_{\mathrm{Sym}(\mathbb{F}_{p})}(H)\leq 
\mathbb{F}_{p}$. Let $h\in H\backslash \mathbb{F}_{p}$. Then $h$ acts on $%
\mathbb{F}_{p}$ as multiplication by a non-identity element, thus $Z_{%
\mathbb{F}_{p}}(h)=\{0\}$, giving $Z_{\mathrm{Sym}(\mathbb{F}_{p})}(H)=1$. $%
\square $

\bigskip

We are ready to prove Theorem \ref{pelda}.

\bigskip

\noindent \textbf{Proof of Theorem \ref{pelda}. }Let $p_{1},p_{2},\ldots
,p_{n}$ be the first $n$ odd primes and let $D=p_{1}p_{2}\cdots p_{n}$. For $%
1\leq i\leq n$ let $m_{i}\in \left\{ 1,\ldots ,D\right\} $ be the (unique)
solution of the congruence system 
\begin{equation*}
m_{i}\equiv \left\{ 
\begin{array}{cc}
1\text{ (}\mathrm{mod}\text{ }p_{j}\text{)} & \text{if }i=j \\ 
2\text{ (}\mathrm{mod}\text{ }p_{j}\text{)} & \text{if }i\neq j
\end{array}
\right. \text{ (}1\leq j\leq n\text{)}
\end{equation*}
and let $S_{i}$ be the set of primes in the arithmetic progression 
\begin{equation*}
\left\{ Dx+m_{i}\mid x\in \mathbb{N}\right\} .
\end{equation*}
Then the $S_{i}$ ($1\leq i\leq n$) are pairwise disjoint.

We claim that $S_{i}$ satisfies both assumptions in Theorem \ref{szamolos} ($%
1\leq i\leq n$). The first assumption follows from the asymptotic form of
Dirichlet's theorem saying $f_{n}(S_{i})=O(n/\log n)$. For the second
assumption let $a$ and $r$ be positive integers; we shall check that the
assumption holds for the arithmetic progression $\left\{ ax+r\mid x\in 
\mathbb{N}\right\} $. Let $p_{j}^{o_{j}}$ be the maximal $p_{j}$-power
dividing $a$ ($1\leq j\leq n$), let 
\begin{equation*}
b=\prod_{j=1}^{n}p_{j}^{o_{j}}
\end{equation*}
and let $a^{\prime }=a/b$ and let $D^{\prime }$ be the least common multiple
of $b$ and $D$. Then $a^{\prime }$ and $D^{\prime }$ are relatively prime,
so there exists a solution $m_{i}^{\prime }$ to the congruence system 
\begin{eqnarray*}
m_{i}^{\prime } &\equiv &m_{i}\text{ (}\mathrm{mod}\text{ }D^{\prime }\text{)%
} \\
m_{i}^{\prime } &\equiv &1\text{ (}\mathrm{mod}\text{ }a^{\prime }\text{)}
\end{eqnarray*}
Since $m_{i}^{\prime }$ and $D^{\prime }a^{\prime }$ are relatively prime,
using Dirichlet's theorem, the set 
\begin{equation*}
S_{i}^{\prime }=\left\{ x\in S_{i}\mid x\equiv m_{i}^{\prime }\text{ (}%
\mathrm{mod}\text{ }D^{\prime }a^{\prime }\text{)}\right\}
\end{equation*}
consists of infinitely many primes. Also there exists $t$ with $%
m_{i}^{\prime }t\equiv r$ ($\mathrm{mod}\text{ }D^{\prime }a^{\prime }$).
Let $s_{1},s_{2},\ldots ,s_{t}$ be distinct elements of $S_{i}^{\prime
}\subseteq S_{i}.$ Then since $a$ divides $D^{\prime }a^{\prime }$, we have 
\begin{equation*}
\sum_{j=1}^{t}s_{j}\equiv m_{i}^{\prime }t\equiv r\text{ (}\mathrm{mod}\text{
}a\text{)}
\end{equation*}
provides the required sum in the second assumption. The claim holds.

Now Theorem \ref{szamolos} implies that there exists a natural number $k$
that can be obtained as a sum of different elements of the $S_{i}$ ($1\leq
i\leq n$). Let $q_{i,j}\in S_{i}$ ($1\leq i\leq n$, $1\leq j\leq l_{i}$) be
different primes satisfying the decompositions 
\begin{equation*}
k=\sum_{j=1}^{l_{i}}q_{i,j}
\end{equation*}

Let $C_{i}$ denote the cyclic group of order $p_{i}$. Let $F_{i,j}=\mathbb{F}%
_{q_{i,j}}$, let 
\begin{equation*}
V_{i}=\bigoplus_{j=1}^{l_{i}}F_{i,j}\text{ and }X_{i}=%
\bigcup_{j=1}^{l_{i}}F_{i,j}
\end{equation*}
Then $p_{i}$ divides $q_{i,j}-1$, so $C_{i}$ embeds into the multiplicative
group of $F_{i,j}$. Let $G_{i}$ be the semidirect product of $V_{i}$ and $%
C_{i}$ acting diagonally on the components $F_{i,j}$. This action defines an
embedding of $G_{i}$ into $\mathrm{Sym}(X_{i})$. Let $G_{i,j}$ denote the
action of $G_{i}$ on $F_{i,j}$. Then $G_{i,j}$ is permutation isomorphic to
a subgroup of $\mathrm{Aff}(p)$ properly containing $\mathbb{F}_{p_{i}}$ and
for $j\neq j^{\prime }$ the permutation groups $G_{i,j}$ and $G_{i,j^{\prime
}}$ are not permutation isomorphic. Applying Lemma \ref{pici}, the
centralizer 
\begin{equation*}
Z_{\mathrm{Sym}(X_{i})}(G_{i})=\bigoplus_{j=1}^{l_{i}}Z_{\mathrm{Sym}%
(F_{i,j})}(G_{i,j})=1
\end{equation*}
is trivial. We showed that the $G_{i}$ ($1\leq i\leq n$) have a permutation
action on $k$ points with trivial centralizer in the full symmetric group $%
\mathrm{Sym}(k)$.

It is easy to see that $G_{i}^{\prime }=V_{i}$ and so $G_{i}/G_{i}^{\prime }$
is cyclic. Trivially, $G_{i}$ is solvable and non-cyclic, so $d(G_{i})=2$.
Also 
\begin{equation*}
\left| G_{i}\right| =p_{i}\prod_{j=1}^{l_{i}}q_{i,j}
\end{equation*}
so for $i\neq i^{\prime }$ the orders of $G_{i}$ and $G_{i^{\prime }}$ are
relatively prime.

We estimate $d\left( \widehat{\Gamma }\right) $ using Lemma \ref{becs} with $%
\mathrm{Sym}(k)$ as target group. We have seen that the $G_{i}$ ($1\leq
i\leq n$) have an embedding into $\mathrm{Sym}(k)$ with trivial centralizer.
Taking into account the trivial permutation representation, this gives 
\begin{equation*}
\left| \mathrm{Hom}(G_{i},\mathrm{Sym}(k))\right| \geq \left| \mathrm{Sym}%
(k)\right| +1
\end{equation*}
which yields 
\begin{equation*}
d(\widehat{\Gamma })\geq \sum_{i=1}^{n}h(G_{i},\mathrm{Sym}(k))\geq n\frac{%
\log (k!+1)}{\log k!}>n
\end{equation*}

The theorem holds. $\square $

\bigskip

\noindent \textbf{Remark. }The upper estimate $d\left( \widehat{\Gamma }%
\right) \leq n+s-1$ does not hold in general, even if we assume that all the 
$G_{i}$ are perfect. Indeed, let the $G_{i}$ ($1\leq i\leq n$) be isomorphic
to $A_{5}$, the alternating group on $5$ letters. Now $\mathrm{Hom}%
(A_{5},A_{5})$ consists of the set of automorphisms and the trivial
homomorphism, so 
\begin{equation*}
h(A_{5},A_{5})=\frac{\log 121}{\log 60}\approx 1.1713
\end{equation*}
implying 
\begin{equation*}
d\left( \widehat{\Gamma }\right) \geq 1.1713n.
\end{equation*}

\end{document}